\documentclass[12pt]{amsart}
\usepackage{a4wide,mathrsfs,graphicx,stmaryrd,pictex, amscd}

\theoremstyle{plain}
\newtheorem{Thm}{Theorem}[section]

\newtheorem{Lem}[Thm]{Lemma}
\newtheorem{Prop}[Thm]{Proposition}
\newtheorem{Cor}[Thm]{Corollary}
\theoremstyle{definition}

\theoremstyle{remark}
\newtheorem{Rmk}[Thm]{Remark}

\theoremstyle{definition}

\newcommand{\R}{\mathbb R}
\newcommand{\Z}{\mathbb Z}

\newcommand{\N}{\mathbb N}
\newcommand{\slr}{\text{SL}_2(\R)}

\title[Cross sections for $\alpha$-CF]{Cross sections for geodesic flows and $\alpha$-continued fractions}
\author{ Pierre Arnoux}
\address{Institut de Math\'ematiques de Luminy (UPR 9016),
         163 Avenue de Luminy, case 907,
         13288 Marseille cedex 09,
         France}
\email{arnoux\@ iml.univ-mrs.fr}
\author{Thomas A. Schmidt}
\address{Oregon State University\\Corvallis, OR 97331, USA}
\email{toms@math.orst.edu}
\date{21 July  2012}

\begin{document}

\begin{abstract}  We adjust Arnoux's coding, in terms of regular continued fractions,  of the geodesic flow on the modular surface to give a cross section on  which the return map is a double cover of the natural extension for the $\alpha$-continued fractions, for each $\alpha \in (0,1]$.     The argument is sufficiently robust to apply to the Rosen continued fractions and their recently introduced $\alpha$-variants.   

\end{abstract}

\maketitle


\section{Introduction}
The $\alpha$-continued fractions introduced by Nakada \cite{Nakada81} 
give a one-dimensional family of interval maps, $T_{\alpha}$  with  $\alpha \in [0,1]$.    This family has been the focus of a recent spate of publications,  focusing especially on the entropy of the maps.    Here we answer a natural question arising in this context:   For $\alpha \in (0,1]$ there is a cross section of the geodesic flow on the modular surface for which $T_{\alpha}$ is a factor.

 Kraaikamp, Schmidt and Steiner \cite{KraaikampSSteiner:12}   determined an explicit  planar model $\mathcal T_{\alpha}: \Omega_{\alpha} \to \Omega_{\alpha}$ of the natural extension of the $\alpha$-continued fraction map for positive $\alpha$.    Denote by  $\mu$ the measure on $\mathbb R^2$ given by $(1 + xy)^{-2} dx dy$.   They also showed that $\mu$ is a $\mathcal T_{\alpha}$-invariant measure on $\Omega_{\alpha}$;   let $h(T_{\alpha})$ denote the {\em entropy} of the interval map with respect to $\mu_\alpha$, the invariant probability measure  on $\Omega_{\alpha}$ induced by $\mu$.  The following is Theorem~2 in \cite{KraaikampSSteiner:12}:  For any  $\alpha \in (0,1]$, 
\begin{equation} \label{e:entropArea}
h(T_{\alpha})\, \mu(\Omega_{\alpha}) = \pi^2/6\,.
\end{equation}
 
 With $ T^1 \mathcal M$ denoting the unit tangent bundle of the modular surface (see below for definitions), we show the following.    

\begin{Thm} \label{t:geoFlow}
For any $\alpha \in (0,1]$,     Equation~\eqref{e:entropArea} holds if and only if there exists
a cross section $\widetilde \Sigma_{\alpha} \subset T^1 \mathcal M$  to the geodesic flow with the following two properties: 
\begin{enumerate}
\item  There is 
 a two-to-one (a.e.) map 
$\pi_{\alpha}: \widetilde \Sigma_{\alpha} \to \Omega_{\alpha}$  that projects the normalized   transverse invariant measure on $\widetilde \Sigma_{\alpha}$  induced by Liouville measure  to $\mu_{\alpha}\,$; 
\item Denote by  $\Phi_{\alpha}$  the first return map to $\widetilde \Sigma_{\alpha} $ of  the geodesic flow, 
the following diagram commutes: 

\[
\begin{CD}
\widetilde \Sigma_{\alpha} @>\Phi_{\alpha}>>\widetilde \Sigma_{\alpha}\\
@V\pi_{\alpha}VV  @VV\pi_{\alpha}V\\
\Omega_{\alpha} @>\mathcal T_{\alpha}>>   \Omega_{\alpha} 
\end{CD}\;.
\]
\end{enumerate}

\end{Thm}

\bigskip
Since a system is by definition a factor of its natural extension,  with Equation~\eqref{e:entropArea}   we have  the following. 
\begin{Cor} \label{t:factor}  For $\alpha \in (0,1]$,   the system defined by $T_{\alpha}$ is a factor of a cross section of the geodesic flow on the modular surface.
\end{Cor}

\bigskip 

 That whenever $T_{\alpha}$ is the factor of a cross section for the geodesic flow there is an equation of the form of  Equation~\eqref{e:entropArea} is well-known to be a consequence of a formula of Abramov, see  Section~\ref{ss:firstReturn} for related ideas.      The main contribution of this paper is thus the construction of the appropriate measurable cross sections.   Even this is simply an adaption of the construction given  in \cite{Arnoux94};  however,   in that setting,  the fact that the  map on the cross section was a {\em first} return was evident.   That this is true in our more general setting requires Equation ~\eqref{e:entropArea}; and, we must admit, careful consideration of the normalizations involved in calculating the entropy of the return map to a cross section.\\

We show that for each $\alpha \neq 0$,   the natural extension of $T_{\alpha}$ is {\em double covered} by a corresponding cross section of the geodesic flow on the modular surface --- it is this cardinality two that accounts for the difference in the constant in Equation~\eqref{e:entropArea} from the volume of the unit tangent bundle of the modular surface.  On the fiber  over an $x$ with $T_{\alpha}(x)$ given by a matrix of determinant $-1$ acting on $x$ in the standard linear fractional manner, the first return of the flow switches sheets,  over $x$ where the map is given by a matrix of determinant $+1$  the flow returns to the same sheet.   The case of $\alpha = 1$ is that of the regular continued fractions, and is thus the case that is treated in  \cite{Arnoux94};  there, above any $x$ the flow changes sheet.   For all other values of $\alpha$, there are $x$ above which one returns to the sheet of departure.   

Our techniques are easily adapted to the special case of $\alpha = 0$.   The interval map $T_0$ is known as the by-excess continued fraction map,  defined on the interval $[-1,0]$ it has an infinite invariant measure given by $(1 + x)^{-1} dx$.  Although one thus cannot define the entropy of $T_0$,    Luzzi and Marmi \cite{LuzziMarmi08} show that its entropy in what they call Krengel's sense is zero.    For all $x$,   $T_0(x)$ is given by some matrix of positive determinant and indeed,  one can extend our discussion here to show that the natural extension of the system of $T_0$ can be expressed as a cross section of flow on the modular surface, with the cross section of infinite measure.\\

The study of the intertwined nature of continued fractions and the geodesic flow on the modular surface has a rich history.    One of the first significant steps was E.~Artin's~\cite{Artin24} coding of geodesics in terms of the regular continued fraction expansions of the real endpoints of their lifts.     Hedlund \cite{Hedlund35}  then used this to show the ergodicity of the geodesic flow on the modular surface.   A few years later,   E.~Hopf proved ergodicity for the geodesic flow on any hyperbolic surface of finite volume, see his reprisal \cite{Hopf:71}.    Taking up on an approach revived by Moeckel \cite{Moeckel82},  Series~\cite{Series85} gave an explicit cross section for the geodesic flow on the modular surface such that the regular continued fraction map is given as a factor.   Adler and Flatto gave a second approach, see especially their  \cite{AdlerFlatto91}.     Quite recently,   
D.~Mayer and co-authors ~\cite{MayerStroemberg08}, ~\cite{MayerMuehlenbruch10}  have given explicit cross sections for the geodesic flow on the surfaces uniformized by the Hecke triangle groups such that a variant, due to Nakada ~\cite{Nakada92},  of the Rosen continued fractions \cite{Rosen54}  is given as a factor.     (The second-named author first learned of such possibilities from A.~Haas, see the related treatment of another variant of the Rosen fractions in \cite{GroechenigHaas96}.)   For much of this history,  further motivation,  and also the work of S.~Katok and co-authors using various means to code geodesics on the modular surface,  see  \cite{KatokUgarcovici07}.  

As stated above, our work here is directly inspired by another approach to identifying cross sections related to continued fractions,  that of Arnoux \cite{Arnoux94}.   That work can be viewed as an elaboration of Veech's \cite{Veech:82} notion of ``zippered rectangles'' in the special case of surfaces of genus one.   Here we  proceed by, in a sense, extracting the algebraic expression for the cross section determined in \cite{Arnoux94}, and adjusting this appropriately so as to pass from the planar models of the natural extensions given by \cite{KraaikampSSteiner:12} to our cross sections.    Thus,  whereas the various cross sections mentioned in the previous paragraph were each constructed in what can be called a geometric fashion,  we use a more algebraic approach.\\

 There has been a great deal of recent interest in the $\alpha$-continued fractions maps,  and especially in the function associating to $\alpha$ the entropy of the corresponding map.  
Nakada  \cite{Nakada81}  computed the entropy of his maps for $\alpha \ge 1/2$.  Kraaikamp \cite{Kraaikamp91} gave a more direct fashion to compute these.   Moussa, Cassa and Marmi \cite{MoussaCassaMarmi:99} gave the entropy for the maps with  $\alpha \in [\sqrt{2}-1, 1/2)$.    Luzzi and Marmi  \cite{LuzziMarmi08}  presented numeric data showing that the entropy function $\alpha \mapsto h(T_{\alpha})$ behaves in a rather complicated fashion as $\alpha$ varies.    Nakada and Natsui ~\cite{NakadaNatsui08} gave explicit intervals on which $\alpha \mapsto h(T_{\alpha})$ is respectively constant, increasing, decreasing.    Carminati and Tiozzo ~\cite{Carminati-Tiozzo:10}  extended this work by describing the intervals involved.     The same authors have revisited the question of the shape of the entropy function in terms of $\alpha$ in \cite{Carminati-Tiozzo:11}.       Luzzi and Marmi \cite{LuzziMarmi08}  also gave strong evidence that the entropy function is continuous. (They furthermore gave numeric evidence that the function is far from trivial,  see the graphs they present in their Figures~4--8).  Tiozzo \cite{Tiozzo:09} showed continuity for $\alpha$ above a certain positive constant, Kraaikamp-Schmidt-Steiner \cite{KraaikampSSteiner:12} show the continuity for positive $\alpha$.    (After that work was completed,  Tiozzo in an updated version of \cite{Tiozzo:09} showed that H\"older continuity holds throughout the full interval.) 

  The proof of  the continuity of $\alpha \mapsto h(T_{\alpha})$ given in \cite{KraaikampSSteiner:12} proceeds by first showing  Equation~\eqref{e:entropArea}  so as then to turn to the key result, the continuity of the function $\alpha \mapsto \mu(\Omega_{\alpha})$  (That the geometry, and even topology, of the $\Omega_{\alpha}$ change drastically is evidenced in Figure~6 of \cite{KraaikampSSteiner:12}.)    The proof of Equation~\eqref{e:entropArea} consists mainly in showing that for any $\alpha >0$,  the intersection $\Omega_{\alpha} \cap \Omega_{1}$ is nontrivial,  and the return maps of each of $\mathcal T_{\alpha}$ and $\mathcal T_1$ to this intersection agree;  from this one uses Abramov's formula and the fact that the equation certainly holds for $\alpha = 1$.       (The proof of continuity of $\alpha \mapsto \mu(\Omega_{\alpha})$ is much harder,  in briefest terms:  the shape of $\Omega_{\alpha}$ is determined mainly by the $T_{\alpha}$-orbits of $\alpha$ and $\alpha-1$, there are intervals of $\alpha$ such that these two orbits meet, and in a common fashion; continuity is fairly straightforward for each of these intervals;  the remaining set of  $\alpha$ is such that the key orbits can still be described symbolically,  the difference in measure for nearby $\alpha, \alpha'$ is bounded above using this description.)

 Luzzi and Marmi  \cite{LuzziMarmi08}  seem to be the first to ask  in print if the $\alpha$-continued fractions arise as a factor of the geodesic flow on the modular surface.      Folklore consensus seems to be that determining whether the $\alpha$-continued fractions are or not factors of a cross section for the geodesic flow remained an important unsolved problem.\\

To show the robustness of our basic argument,  in Section~\ref{s:RosenCF} we sketch a proof that the Rosen continued fraction maps and certain of their $\alpha$-type variants are also factors of cross sections of geodesic flow over the corresponding hyperbolic surface.   Rosen \cite{Rosen54} introduced his analogs of the nearest integer continued fractions to assist in the study of the Hecke triangle Fuchsian groups, see below for definitions.     There has been increasing interest in the past few decades in these continued fractions,   see \cite{SchmidtSheingorn95} for a sketch of the history to the mid-1990s,  and the aforementioned  \cite{MayerStroemberg08}, \cite{MayerMuehlenbruch10} for some of the more recent history;   the $\alpha$-Rosen continued fractions, that we treat as well,  were introduced by Dajani, Kraaikamp and Steiner \cite{DajaniKraaikampSteiner09} in 2009. \\

In light of the extremely complicated geometry and topology of the 
various cross sections that we determine, it would be interesting to know if the more geometric approaches to the construction of cross sections giving interval maps as factors could be successful in our setting.

\bigskip 
\noindent
{\bf Thanks} The second-named author takes great pleasure in thanking Cor Kraaikamp and Wolfgang Steiner for introducing him to the intricacies of the $\alpha$-continued fractions and for comments on this continuation of joint work.    Thanks also go to Corinna  Ulcigrai for mentioning the interest in the question of relating cross sections and the $T_{\alpha}$,  and to an anonymous referee who insisted on the importance of it.

\section{Background} 

\subsection{$\alpha$-continued fractions} 
For $\alpha \in [0,1]$, we let $\mathbb{I}_{\alpha} : = [\alpha-1,
\alpha]$ and define the map $T_{\alpha}:\, \mathbb{I}_{\alpha} \to \mathbb{I}_{\alpha}$ by
\[
T_{\alpha}(x) := \left| \frac{1}{x} \right| -  \left \lfloor\,  \left| \frac{1}{x} \right| + 1 -\alpha \right\rfloor,\ \mbox{for}\ x \neq 0\,; \ T_{\alpha}(0) :=0 \,.
\]
For $x \in \mathbb{I}_{\alpha}$, put
\[
\varepsilon(x) := \left\{\begin{array}{cl}1 & \mbox{if}\ x \ge 0\,, \\ -1 & \mbox{if}\ x<0\,,\end{array}\right. \quad \mbox{and} \quad d_{\alpha}(x) := \left\lfloor \left| \frac{1}{x} \right| + 1 - \alpha \right\rfloor\,,
\]
with $d_\alpha(0) = \infty$.

Furthermore, for $n \geq 1$, put
\[
\varepsilon_n = \varepsilon_{\alpha,n}(x) := \varepsilon(T^{n-1}_\alpha (x)) \quad \mbox{and}\quad d_n = d_{\alpha,n}(x) := d_\alpha(T^{n-1}_\alpha (x)).
\]
This yields the \emph{$\alpha$-continued fraction}  expansion of
$x\in\R$\,:
\[
x = d_0+ \dfrac{\varepsilon_1}{d_1 + \dfrac{\varepsilon_2}{d_2 +
\cdots}} 
\,,
\]
where $d_0\in\Z$ is such that $x-d_0\in \mathbb{I}_{\alpha}$.
(Standard convergence arguments justify equality of $x$ and its
expansion.) These include the regular continued fractions, given by
$\alpha = 1$ and the nearest integer continued fractions, given by
$\alpha = 1/2$.

The standard number theoretic planar map associated to continued fractions is defined by
\[
\mathcal{T}_\alpha(x,y) := \bigg(T_\alpha(x), \frac{1}{d_{\alpha}(x)+ \varepsilon(x)\,y}\bigg)\, \quad (x \in \mathbb{I}_\alpha,\ y \in [0,1])\,,
\]
and $\Omega_{\alpha}$ is the closure of the orbits of $(x,0)$ with $x \in \mathbb{I}_\alpha$.     
 
 Recall that $\mu$ is given by $(1 + xy)^{-2} dx dy$ and that $\Omega_{\alpha}$ is the planar region determined by \cite{KraaikampSSteiner:12}.  
Let $\mu_\alpha$~be the probability measure given by normalizing $\mu$  on~$\Omega_\alpha$, and $\nu_\alpha$~the  marginal probability measure obtained by integrating $\mu_\alpha$ over the fibers~$\{x\} \times \{y \mid (x,y) \in \Omega_\alpha\}$, $\mathscr{B}_\alpha$~the Borel $\sigma$-algebra of~$\mathbb{I}_\alpha$, and $\mathscr{B}_\alpha'$ the Borel $\sigma$-algebra of~$\Omega_\alpha$.
That  $(\Omega_\alpha, \mathcal{T}_\alpha, \mathscr{B}_\alpha', \mu_\alpha)$ is a natural extension of $(\mathbb{I}_\alpha, T_\alpha, \mathscr{B}_\alpha, \nu_\alpha)$ is shown in \cite{KraaikampSSteiner:12}.

\subsection{Geodesic flow on the modular surface}

Much of the following can be found in Manning's chapter \cite{Manning91} in the text \cite{BedfordKeaneSeries91}.\\  

    Using the M\"obius action of $\text{SL}_2(\mathbb R)$ on the Poincar\'e upper-half plane $\mathbb H$,  by identifying a matrix with the image of $z=i$ under it,   we can identify $\text{SL}_2(\mathbb R)/ \text{SO}_2(\mathbb R)$ with $\mathbb H$.   Similarly,   $\text{PSL}_2(\mathbb R) = \text{SL}_2(\mathbb R)/\pm I$ can be identified with the {\em unit tangent bundle} of $\mathbb H$.   We are most interested in the {\em modular surface},  $\mathcal M =  \text{SL}_2(\mathbb Z)\backslash  \mathbb H$.     The unit tangent bundle of the modular surface, $T^1 \mathcal M$,   can be identified as 
\[T^1 \mathcal M =  \text{PSL}_2(\mathbb Z)\backslash \text{PSL}_2(\mathbb R)\,.\]

The geodesic flow in our elementary setting is a map on a surface's unit tangent bundle:  Given a time $t$ and a unit tangent vector $v$,   since the unit tangent vector uniquely determines a geodesic passing through the vector's base point,  we can follow this geodesic for arclength t in the direction of $v$,   the unit vector that is tangent to the geodesic at the end point of the geodesic arc is the image,  $g_t(v)$, under the geodesic flow.    The hyperbolic metric on $\mathbb H$ corresponds to an element of arclength satisfying $ds^2 = (dx^2 + dy^2)/y^2$ with coordinates $z = x + i y$.   In particular,  for $t>0$,   the points $z = i$ and $w = e^t i$ are at distance $t$ apart.  Since $\text{SL}_2(\mathbb R)$ acts by  isometries on $\mathbb H$,   the geodesic flow on its unit tangent bundle is given by sending $A\in \text{PSL}_2(\mathbb R)$ to $A g_t$,   where $g_t = \begin{pmatrix} e^{t/2}&0\\0&e^{-t/2}\end{pmatrix}$.     Similarly,  on $T^1 \mathcal M$   one sends the class represented by $A$ to that represented by $A g_t$.   

There is a natural measure on the unit tangent bundle $T^1 \mathbb H$:   {\em Liouville measure} is given  as the product of the hyperbolic area measure on $\mathbb H$ with the length measure on the circle of unit vectors at any point.   This measure is $\text{SL}_2(\mathbb R)$-invariant,  and induces a finite measure on $T^1 \mathcal M$.   With the standard choice of this Liouville measure, the volume of $T^1 \mathcal M$ is  $\text{vol}(T^1 \mathcal M) = \pi^2/3$, see say  \cite{Arnoux94}. {\em Normalized Liouville measure} is the corresponding probability measure  on  $T^1 \mathcal M$.  Gurevich and Katok \cite{GurevichKatok:01} cite a fairly general result of Sullivan \cite{Sullivan:84}  to state that the entropy with respect to  {\em normalized} Liouville measure of the geodesic flow  on the modular surface  is one.\\

A measurable {\em cross section} for the geodesic flow is a subset of the unit tangent bundle through which almost every geodesic passes tranversely and infinitely often.   The {\em first return map} to a cross section $\mathcal C$  is the self-map $\Phi: \mathcal C \to \mathcal C$ given by    $\Phi(v) = g_{\tau(v)}(v)$ with $\tau(v)>0$ minimal such that this image is indeed in $\mathcal C$.   By definition of $\mathcal C$,  this map is defined almost everywhere on $\mathcal C$.     By a celebrated result of Ambrose \cite{Ambrose41},   since the geodesic flow leaves normalized Liouville measure invariant,  there is a measure $\lambda$ on $\mathcal C$ invariant under the the first return map such that normalized Liouville measure locally factors as the product of $\lambda$ with Lebesgue measure along geodesics.  In subsection~\ref{ss:mainResults},  we will refer to this phenomenon as an {\em Ambrose factorization}.

Due to a result of Abramov \cite{Abramov:59b},
one defines the entropy of a flow $\{\phi_t\,\vert\, t \in \mathbb R\}$ on a probability space as the entropy of the time $t=1$ map, $h(\phi_1)$.   
A celebrated formula of Abramov  \cite{Abramov:59a} then gives that the {\em entropy of the  first return map} $\Phi$ to a cross section $\mathcal C$ is 
\[ h(\Phi) = h( \phi_{1}) \int_{\mathcal C}\, \tau(v)\, d \lambda_1 = \dfrac{h( \phi_{1})}{\lambda(\mathcal C)}\,  \int_{\mathcal C}\, \tau(v)\, d \lambda\,,\]
where $\lambda_1$ is the probability measure on $\mathcal C$ induced by $\lambda$.   

The integral of the return times over a cross section gives the volume of the total space.   In our setting,  this volume is one since we considered normalized Liouville measure;  Sullivan's result gives that $h( \phi_{1}) = 1$.   Thus,  we have 
\begin{equation}\label{e:entropyReturn}
h(\Phi) = 1/ \lambda(\mathcal C)\,.
\end{equation}

\section{Louville measure as a Lebesgue measure} \label{ss:arnouxHaar}        
Recall that any locally compact Hausdorff topological group has an invariant measure for left multiplication,  this {\em Haar measure} is unique up to scaling.   
In this section,  we derive explicit formulas for the Haar measure on $\text{SL}_2(\mathbb R)$.   These formulas are presumably well known,  and can verified by specializing known expressions for the Haar measure on various classes of Lie groups.  Here we derive them in a straightforward, elementary fashion.   This derivation then allows us to  identify a transversal to the geodesic flow --- see the form of the matrices given  in Lemma~\ref{l:towardsTransv} --- that is key to our construction of cross sections.

\smallskip 
\noindent
{\bf Notation.} 
Throughout this section,  in order to use standard naming conventions of entries of matrices,    $\alpha$ denotes any real number. 

 \bigskip 
Let $G_{\gamma}  \subset \text{SL}_2(\mathbb R)$ be  the set
\[G_{\gamma} = \bigg\{ \begin{pmatrix} \alpha & \beta \\ \gamma & \delta  \end{pmatrix}\,\vert\, \gamma \neq 0\,\bigg\}\,.\]

\begin{Prop}\label{p:haar} The Haar measure $h$ on $G_{\gamma} $ is given, up to a constant, by 
\[d h = \frac{d\alpha\, d\gamma \, d\delta\,}{|\gamma|}\,\]
\end{Prop}

\begin{Lem}   The natural left action of the group $\slr$ on the set $\mathcal M_2(\R)$ of $2\times 2$ real matrices preserves the measure induced by Lebesgue measure on $\R^4$ under the standard identification of $\mathcal M_2(\R)$ with $\R^4$.
\end{Lem} 

\begin{proof} Left multiplication on $\mathcal M_2(\R)$ by a matrix  $M=\begin{pmatrix} a&b\\c&d \end{pmatrix}\in \mathcal M_2(\R)$,   defines a linear transformation,  whose matrix with respect to the canonical basis is
\[ \begin{pmatrix} a&0&b&0\\0&a&0&b\\ c&0&d&0\\0&c&0&d \end{pmatrix}\,.\]
The determinant of this new matrix is  $(ad-bc)^2$, and thus equals $1$ when  $M\in\slr$. That is,  the left multiplication preserves 
Lebesgue measure.
\end{proof}

\begin{proof}[Proof of Proposition~\ref{p:haar}]
Let   $\Delta$ be the  determinant $\alpha\delta-\beta\gamma$. 
On the subset of matrices with $\gamma \neq 0$,  the map 
$(\alpha,\Delta, \gamma, \delta)\mapsto 
(\alpha,(\alpha \delta - \Delta)/\gamma,\gamma, \delta\,)$ 
is injective,  and thus $(\alpha,\gamma,\delta, \Delta)$ gives local coordinates on this subset.    The Jacobian of this map equals $|\gamma^{-1}|$, and thus Lebesgue measure  $m$ on 
$\mathcal M_2(\R)\equiv \R^4$ is given in these coordinates by 
\[dm=d\alpha\,d\beta\,d\gamma\,d\delta=
\frac{d\alpha\, d\gamma\,d\delta\, d\Delta}{|\gamma|}\,.\]

Now, left multiplication by elements of $\slr$ preserves $\Delta$ on 
$\mathcal M_2(\R)$.    Therefore,   restricting to the hypersurface where $\Delta=1$ (that is to $\slr$),  we find that the measure 
$\frac{d\alpha\,d\gamma\, d\delta}{| \gamma |}$ is invariant by left multiplication.  This of course implies that the measure is a Haar measure. 
\end{proof}

\begin{Rmk}   There are obviously three related expressions:  
$\frac{d\alpha\,d\beta\,d\delta}{|\beta |}$ 
etc.  By a Jacobian change-of-coordinates calculation, one easily shows  that each of these expressions defines the same measure on the intersection of their respective domains of definition. 

We emphasize that the above  directly proves the invariance of the measure under left (and, for that matter, also right) multiplication.   That is,  we have naively solved for an explicit expression of Haar measure.
\end{Rmk}

\begin{Lem}\label{l:towardsTransv}   Let  $G^{+}_{\gamma}$  be the connected component of $G_{\gamma}$ defined by $\gamma >0$, and let 
$\Sigma \subset G^{+}_{\gamma}$ be defined by $\gamma = 1$.  Consider local coordinates $x, y$ by letting each 
$A \in \Sigma$ be given as   
\[ A =  \begin{pmatrix} x & x y-1 \\ 1 & y  \end{pmatrix}\,.\]
Then $G^{+}_{\gamma}$ has local coordinates $(x,y,t)$ by way of 
\[ M = A \, g_t\,\]
with $A \in \Sigma$ and $g_t$ as above.    Furthermore, $dx\, dy\, dt$ gives  Haar measure restricted to  $G^{+}_{\gamma}$.
 
\end{Lem} 

\begin{proof}   Suppose that $M=\begin{pmatrix} \alpha & \beta \\ \gamma & \delta  \end{pmatrix} \in G^{+}_{\gamma}$.    Then letting $t = 2 \log \gamma$ gives  $A = M\, g_{-t} \in \Sigma$.      Clearly,  the set $A \, g_t$ for all $A\in \Sigma$,  $t \in \R$ comprises all of $G^{+}_{\gamma}$.      

The map $(x,y,t) \mapsto (\alpha, \gamma, \delta) = (x e^{t/2}, e^{t/2}, y e^{-t/2}\,)$ has Jacobian of absolute value $\gamma/2$.  Since Haar measure restricted to $G^{+}_{\gamma}$ is (any nonzero constant times)  $d h=\frac{d\alpha\, d\gamma \, d\delta\,}{\gamma}$, we deduce that 
$dx\, dy\, dt$ does give Haar measure here.  (The change in normalization constant by a factor of $2$ is innocuous here.)
\end{proof}

For an element $M \in \slr$,  let $[M]$ denote the corresponding element of $\text{PSL}_2(\R)$.  Then since Haar measure $\slr$ induces Liouville measure on $T^1 \mathbb H$,    and $G^{+}_{\gamma}$ projects one-to-one to a set of full measure in $\text{PSL}_2(\R)$, we have the following. 

\begin{Thm}\label{t:liouville} Under the identification of $T^1 \mathbb H$ with 
$\emph{PSL}_2(\mathbb R)$,   Liouville measure on the full measure set $\{\, [M]\,|\, M \in G_{\gamma}\,\}$ is proportional to 
$dx\, dy\, dt$ where $x,y,t$ are as above. 
\end{Thm}
  
(We identify the above implied constant of proportionality in  Subsection~\ref{ss:firstReturn}.)

\bigskip
One can give another parametrization by remarking that a matrix $\begin{pmatrix} \alpha & \beta \\ \gamma & \delta  \end{pmatrix}$ in $\text{PSL}_2(\mathbb R)$ corresponds to an horizontal tangent vector pointing to the right if and only if $\gamma=\delta$, and to the left if and only if $\gamma=-\delta$. The set of such matrices gives a section, since any geodesic which is not vertical has a highest point, where the tangent vector is horizontal. One can parametrize such a geodesic by its origin $\frac {\beta}{\delta}$ and its extremity $\frac{\alpha}{\gamma}$. A simple computation shows that the set of such matrices, in the case $\gamma=\delta$, is parametrized by $\begin{pmatrix}\frac{X}{\sqrt{X-Y}} & \frac{Y}{\sqrt{X-Y} }\\ \frac{1}{\sqrt{X-Y}} &\frac{1}{\sqrt{X-Y}} \end{pmatrix}$

Applying the geodesic flow  as above gives a parametrization of half of the group of the form 

\[\begin{pmatrix}\frac{Xe^{\frac t2}}{\sqrt{X-Y}} & \frac{Ye^{\frac {-t}2}}{\sqrt{X-Y} }\\ \frac{e^{\frac t2}}{\sqrt{X-Y}} &\frac{e^{\frac {-t}2}}{\sqrt{X-Y}} \end{pmatrix}\,.
\]

A Jacobian computation proves that the Haar measure, in these coordinates, is proportional to $\frac {dX\, dY\, dt}{(X-Y)^2}$. Note that the measure   $\frac {dX\, dY}{(X-Y)^2}$ has a clear geometric signification : the base point of the corresponding unit tangent vector is $\frac {X+Y}2+i \frac{X-Y}2$, so this measure is proportional to the measure defined by the hyperbolic metric on the upper half-plane.

A small modification consists in changing $Y$ into $-1/Y$, giving the following result. 

\begin{Prop} Any element of $\emph{PSL}_2(\mathbb R)$ such that $\gamma$ and $\delta$ have same sign can be written in a unique way (up to sign):
\[\begin{pmatrix}\frac{Xe^{\frac t2}}{\sqrt{X-Y}} & \frac{e^{\frac {-t}2}}{-Y\sqrt{X-Y} }\\ \frac{e^{\frac t2}}{\sqrt{X-Y}} &\frac{e^{\frac {-t}2}}{\sqrt{X-Y}} \end{pmatrix}\,.\]

In these coordinates, the Haar measure can be written 
\[\frac {dX\, dY\, dt}{(1+XY)^2}\,.\]
\end{Prop}

Thus, we have sketched a geometric interpretation of the density function for the measure $\mu$ on $\mathbb R^2$ that is sometimes referred to as the ``standard  number theoretic'' measure.   (This latter name comes from the fact  that the orbits of points of the form $(x,0)$ in planar models of natural extensions with $\mu$ an invariant measure reveal Diophantine approximation properties of the related interval maps.)  Note that a similar derivation of an invariant measure for Rosen continued fraction from Liouville measure for the geodesic flow has been given by Mayer and Stroemberg in the section 5 of \cite{MayerStroemberg08}.

It is easy to give explicit formulas to change between these systems of coordinates; we use such a change of coordinates in the next section, the function $\mathcal Z$, to transform the invariant measure to Lebesgue measure.

\section{The cross section}   The goal of this section is, for each $\alpha \in (0,1]$, to define a subset $\widetilde \Sigma_{\alpha}$ of the unit tangent bundle of the modular surface and to give an explicit self-map induced by the geodesic flow. 

\subsection{An alternate natural extension}\label{ss:altNatExt} 
Let the function $\mathcal Z: \mathbb R^2 \setminus \{(x,y)\,\vert\, y = -1/x\}  \to \mathbb R^2$ be defined by $\mathcal Z(x,y) = (x, y/(1+ x y)\,)$ and let  $\Sigma_{\alpha} = \mathcal Z(\Omega_{\alpha})$.        An elementary calculation shows that $\mathcal Z$ conjugates $\mathcal T_{\alpha}$ to the bijection:

\begin{align} \label{e:biject}
\Sigma_{\alpha} &\to \Sigma_{\alpha}\\
(x,y) &\mapsto  (\, f_{\alpha}(x),  \, \varepsilon(x) \; x (1-x y)\,)\,.\notag
\end{align}

\noindent
A second elementary calculation shows that $\mu$ projects to give Lebesgue measure, which is thus invariant for this induced map.  (Of course, one can just as easily directly verify this invariance.)

\subsection{Definition and main results} \label{ss:mainResults}

Given   $(x, y) \in \mathbb R^2$,  let

\begin{align}\label{e:basicMatShapes}  A_{-1}(x,y) &= \begin{pmatrix} 1&y\\-x & 1 - x y   \end{pmatrix}\;\;\mbox{and} \notag \\
\\
A_{+1}(x,y) &= \begin{pmatrix} x & 1 - x y \\ -1 & y  \end{pmatrix}\,.
\notag
\end{align}

Define the following subsets of $\text{PSL}_2(\mathbb R)$:   
\[ \mathcal A_{\alpha, -1} = \{\;\left[ A_{-1}(x,y)\right] \,     \vert \, (x, y) \in \Sigma_{\alpha}\}\;\mbox{and} \;  \mathcal A_{\alpha, +1} = \{ \; \left[ A_{+1}(x,y)\right] \, \vert \, (x, y) \in \Sigma_{\alpha}\}\,.   \]

In Lemma~\ref{l:distinctClasses} we show that almost all  $A, A'$ in $\bigcup_{\sigma \in \{-1,+1\} }\; \mathcal A_{\alpha, \sigma}$  lie in distinct  $ \text{PSL}_2(\mathbb Z)$-orbits.  We let $\widetilde \Sigma_{\alpha, \sigma}$   denote the set of the classes in 
$ \text{PSL}_2(\mathbb Z)\backslash \text{PSL}_2(\mathbb R)$ represented by the various elements of the $\mathcal A_{\alpha, \sigma}$ and let 
\[ \widetilde \Sigma_{\alpha} = \bigcup_{\sigma \in \{-1,+1\} }\, \widetilde \Sigma_{\alpha, \sigma}\,.\]
 Let $\ell$ on  $\widetilde \Sigma_{\alpha}$ be induced by Lebesgue measure  on $\Sigma_{\alpha}$.    The previous section applies to show that the geodesic flow from $\widetilde \Sigma_{\alpha}$ gives an Ambrose factorization of (some positive constant multiple of) Liouville measure of the form $\ell \times dt$.

Denote a point of $\widetilde \Sigma_{\alpha}$ by $(x, y, \sigma)$ in the obvious fashion.   Our main result is given by the following two statements; we prove Theorem~\ref{t:First}  in Section~\ref{s:firstReturn}.

\begin{Thm} \label{t:First}
For any $\alpha \in (0,1]$,  the first return map of the geodesic flow to $\widetilde \Sigma_{\alpha}$ is given by 
\[
\begin{aligned}  \Phi_{\alpha}: \;\;\; \;\;\;\;\widetilde \Sigma_{\alpha}\;\;\; &\to \;\;\;\;\widetilde \Sigma_{\alpha}\\
(x, y, \sigma) &\mapsto  (\, f_{\alpha}(x),  \,\varepsilon(x) \;  x (1-x y),  - \varepsilon(x) \sigma\,)\;.
\end{aligned} 
\]
\end{Thm}

By the discussion in Subsection~\ref{ss:altNatExt},  the following holds.
\begin{Lem} \label{l:cover}
For  $\alpha \in (0,1]$, let
\[
\begin{aligned}
\pi_{\alpha}: \;\;\; \;\;\;\;\widetilde \Sigma_{\alpha}\;\;\; &\to \;\;\;\Omega_{\alpha}\\
(x, y, \sigma) &\mapsto (\,x, \, y/(1+x y)\,)\;.
\end{aligned}
\]

   Then $\pi_{\alpha}$ is a $2:1$ surjection,  and the following is a commutative diagram.

\[
\begin{CD}
\widetilde \Sigma_{\alpha} @>\Phi_{\alpha}>>\widetilde \Sigma_{\alpha}\\
@V\pi_{\alpha}VV  @VV\pi_{\alpha}V\\
\Omega_{\alpha} @>\mathcal T_{\alpha}>>   \Omega_{\alpha} 
\end{CD}
\]
\end{Lem}

\section{ First return map}\label{s:firstReturn}

In this section we give the remaining steps to prove the main result, in particular showing that $\Phi_{\alpha}$ is a {\em first} return map.  To show that this step is in general necessary,   in the final subsection  we give an example (which can easily be generalized) of an interval map that is not given as the factor of first return to a cross section for the geodesic flow.

\smallskip 
\noindent
{\bf Notation.} 
We use asterisks to denote entries of a determinant one matrix that are not germane to the argument at hand. 

 \bigskip

\subsection{Classes are distinct} 
We first show that the classes parametrized by the various  $(x, y, \sigma) \in \widetilde \Sigma_{\alpha}$ are almost always distinct.     This implies that
\begin{equation}\label{e:twiceMu}
 \ell(\widetilde \Sigma_{\alpha}) = 2 \mu(\Omega_{\alpha})\,.
 \end{equation}

\bigskip
\begin{Lem} \label{l:distinctClasses}
Let $\alpha \in (0,1)$.   For almost all distinct  $A, A'$ in $\bigcup_{\sigma = \pm 1}\; \mathcal A_{\alpha, \sigma}$ the $ \emph{PSL}_2(\mathbb Z)\backslash \emph{PSL}_2(\mathbb R)$ classes represented by $A, A'$ are distinct.  
\end{Lem}

 \begin{proof}  {\bf case 1.}\quad   Suppose that both $A, A' \in \mathcal A_{\alpha, -1}$.   
 
    If  $M\in \text{SL}_2(\mathbb Z)$ is such that 
$M A = \pm A'$,  then  writing  $A, A'$ in accordance with  \eqref{e:basicMatShapes} and multiplying on the right by $A^{-1}$ gives 
 \[ M = \pm  \begin{pmatrix} (1-x y) + xy'&- y +  y'\\ *&*\end{pmatrix}\,.
 \]
Since $M \in \text{SL}_2(\mathbb Z)$, we find that $y' - y \in \mathbb Z$.   Substitution into the $(1,1)$-element easily leads to the conclusion that $x\in \mathbb Q$.

\bigskip 
\noindent
{\bf case 2.}\quad  If both $A, A' \in \mathcal A_{\alpha, +1}$, then we find 
 \[ M = \pm  \begin{pmatrix} *&*\\ -y+ y'&1-x y + xy'\end{pmatrix}\,,
 \]
 and again $x \in \mathbb Q$.

\bigskip 
\noindent
{\bf case 3.}\quad   If $A \in \mathcal A_{\alpha, -1}$ and $A' \in \mathcal A_{\alpha, +1}$, then again $x \in\mathbb Q$,  as 
 \[ A' A^{-1}  = \pm  \begin{pmatrix} *&*\\  1+ x y + xy'&y+ y'  \end{pmatrix}\,.
 \]
 
 Thus,  in all cases, symmetry shows that distinctness holds unless both $x, x'$ belong to the countable set  $\mathbb Q$.
\end{proof}

\subsection{Flow does give $\Phi_{\alpha}$ as in Theorem~\ref{t:First}}

\begin{Lem} \label{l:flowFormula}
For each $A =  (x, y,  \sigma) \in \widetilde \Sigma_{\alpha}$, the geodesic flow for time $t = -2 \log |x|$ sends $A$  to $(\, f_{\alpha}(x),  \, \varepsilon(x) \; x (1-x y),   - \varepsilon(x) \sigma\,)$.
\end{Lem}

\begin{proof}  Consider 
$A = \begin{pmatrix} 1&y\\-x & 1 - x y  \end{pmatrix}$ representing a class of $\widetilde \Sigma_{\alpha, -1}$.  Since $t = -2 \log |x|$ gives  $e^{t/2} = \varepsilon/x$ where $\varepsilon = \varepsilon(x)$, the geodesic flow  gives $A g_t = \begin{pmatrix} \varepsilon/x&\varepsilon x y\\-\varepsilon & \varepsilon x(1 - x y)   \end{pmatrix}$.   When $\varepsilon = -1$,  this is equivalent to 
\[\begin{pmatrix} 0&1\\ -1& d   \end{pmatrix} A g_t = \begin{pmatrix} 1& - x(1 - x y) \\-(-d - 1/x)&* \end{pmatrix}\,,\]
where we choose $d = d_{\alpha}(x)$.    When $\varepsilon = 1$,  we choose a different $\text{PSL}(2, \mathbb Z)$-orbit representative of  $Ag_t$, to wit: 
\[\begin{pmatrix} 1&d\\ 0& 1   \end{pmatrix} A g_t = \begin{pmatrix} -d + 1/x&*\\-1&  x(1 - x y)   \end{pmatrix}\,.\]

Similarly,  for $A = \begin{pmatrix} x &  1 - x y   \\ -1 & y  \end{pmatrix} $ representing a class in  $\widetilde \Sigma_{\alpha, +1}$, if $\varepsilon = -1$,  we find 

\[ \begin{pmatrix}  d&-1\\ 1&0\end{pmatrix} A g_t = \begin{pmatrix}   -d-1/x &*\\-1 &-x(1-xy) \end{pmatrix}\,.
\]
If $\varepsilon = 1$,  we have 
\[ \begin{pmatrix}  1&0\\ d&1\end{pmatrix} A g_t = \begin{pmatrix}  1&x(1-xy)\\-(-d + 1/x) &* \end{pmatrix}\,.
\]

  Due to the bijection of \eqref{e:biject}, and the definition of the $\widetilde \Sigma_{\alpha, \sigma}$,  the result holds. 
\end{proof} 

\subsection{The map $\Phi_{\alpha}$ is given by the {\em first} return}\label{ss:firstReturn}
Let $\ell_{\alpha}$ be the probability measure on $\widetilde \Sigma_{\alpha}$ induced by Lebesgue measure $\ell$.   
The system  $(\widetilde \Sigma_{\alpha}, \ell_{\alpha}, \Phi_{\alpha})$ is a skew product with finite fiber over $(\Omega_{\alpha}, \mu_{\alpha}, \mathcal T_{\alpha})$.   Their entropies are thus equal:   $h(\Phi_{\alpha}) = h( \mathcal T_{\alpha})$.   Since a system and its natural extension have the same entropy,  Equation~\eqref{e:entropArea}  gives $h(\Phi_{\alpha}) \mu(\Omega_{\alpha}) = \pi^2/6$.   

By Equation~\eqref{e:entropyReturn}, Theorem~\ref{t:liouville} and Equation~\eqref{e:twiceMu},  the first return map of the geodesic flow to $\widetilde \Sigma_{\alpha}$ has entropy equal to  
\[\dfrac{1}{ \lambda(\widetilde \Sigma_{\alpha})} = \dfrac{1}{ \kappa\, \ell(\widetilde \Sigma_{\alpha})} = \dfrac{1}{ 2 \kappa \,\mu( \Omega_{\alpha})}\,,\]
where $\kappa$ is a  normalizing constant  independent of $\alpha$.   But, $(\widetilde \Sigma_{1}, \ell_{1}, \Phi_{1})$ {\em is} given by the first return of the geodesic flow (by, say, \cite{Arnoux94}), and we thus find that     $1/\kappa= \pi^2/3$.   Therefore, for all positive $\alpha$, we find that the entropy of the first return map to $\widetilde \Sigma_{\alpha}$ equals $h(\Phi_{\alpha})$.     But,  $\Phi_{\alpha}$ is given (locally) by powers of this first return map,  and thus  the equality of $h(\Phi_{\alpha})$  with the entropy of the first return map holds   if and only if the two maps are equal a.e.   --- see \S 10.6 Theorem 2 (3) of \cite{CornfeldFominSinai82}.    That is,  $\Phi_{\alpha}$ is the {\em first} return map.

\subsection{Non-first return interval maps }  To reassure the reader that our efforts in the previous subsection are not absurd,  we sketch the existence of an interval map that is piecewise fractional linear with integral coefficients  with the planar model of its natural extension double covered by a cross section for the geodesic flow on $T^1 \mathcal M$ but for which the map is not given by the $n^{\text{th}}$ return map for any single $n$.   

We begin with the regular continued fractions.    For simplicity,  let $I = \mathbb{I}_1$ and  $T = T_1$.     Denote cylinder sets in the usual manner:   $\Delta[\,a_1, \dots , a_n\,]$ is the subset of the unit interval of elements whose first $n$ partial quotients are the $a_i\,$.  In particular,  the $T$-image of    $\Delta[\,a_1, \dots , a_n\,]$ is $\Delta[\,a_2, \dots , a_n\,]$ whenever $n \ge 2$.  
Using the  partition of $I\,$  given by
\[   A_1 = \; \Delta[1]\,, \;\mbox{and the various}\;\;  A_n = \;\bigcup_{k>1, l>1}\; \Delta[k,\underbrace{1, \dots, 1}_{n-2 \,  \text{times}}, l]\;\; \mbox{with}\; n\ge 2 \;,\]

\medskip 
\noindent
let 
\[ g : I \to I\] 
be defined by 
\[ g(\, A_n\,) =  T^n(A_n)\,.\] 
That is, on each $A_n$ we define $g$ to be the $n$-fold composition of $T$ with itself.

\medskip  
Again simplifying notation, let $\mathcal T: \Omega \to \Omega$ be the planar natural extension for $T\,$.  For each $n \ge 1$, define  $\mathcal A_n$ to be the subset of $\Omega$  lying over $A_n\,$.    Let  $\mathcal G: \Omega \to \Omega$ be defined by letting $\mathcal G$ restricted to $\mathcal A_n$ be the $n$-fold composition of $\mathcal T$ with itself. 

As usual, for simplicity's sake, we refer to dynamical systems merely by space and function,  each time the  mentioned have sigma-algebra of Borel subsets;  invariant measures in the following are (the normalization of) $d \mu = (1 + x y)^{-2} \, dx dy$ and its marginal measure.   We let $\widetilde \Sigma =   \widetilde \Sigma_1$ as defined in Subsection~\ref{ss:mainResults}.

\begin{Lem}   The map $\mathcal G: \Omega \to \Omega$  is a natural extension for $g$.   Furthermore,  $ \widetilde\Sigma \subset \mathcal T^1 \mathcal M$ is partitioned by sets indexed by $\mathbb N$ such that the map on $ \widetilde\Sigma$, sending any element of the $n$th partition set to its $n$th return under the geodesic flow to $\widetilde\Sigma$, gives a double cover of $\mathcal G$.   
\end{Lem}
 
\begin{proof}  The union of the $\mathcal A_n$ is clearly all of $\Omega\,$,  up to a set of measure zero.    We claim that also the union of the  $\mathcal T^{n-1}(\mathcal A_{n}\,)$ is $\Omega\,$, up to measure zero.      To see this, first note that for $n \ge 2$,  
\[T^{n-1} A_n =  I \setminus \Delta[1] = [0, 1/2)\,.\]
Now,  let $N_d = \begin{pmatrix}  0&1\\1&d\end{pmatrix}$; thus,   for $(x,y) \in \mathcal A_n$ with $x \in \Delta[k]$,   
\[ \mathcal T^{n-1} (x,y) = (\, T^{n-1}(x), \,N_{1}^{n-2}N_k\cdot y\,)\,.\]
A proof by induction shows that for $j \ge 1$, 
\[ N_{1}^{j} = \begin{pmatrix}  f_{j-2}&f_{j-1}\\f_{j-1}&f_{j}\end{pmatrix}\,\]
where $f_0 = 0, f_1 = 1, f_2 = 1, f_{j+2} = f_{j+1}+f_j$ is the Fibonacci sequence.    Since $N_k\cdot y = 1/(y+k)$, we have that 
\[\mathcal T^{n-1} \mathcal A_n = [0,1/2) \times N_1^{n-2}\cdot[0,1/2)\,.\]
     But, $N_{1}^{n-2}\cdot 0 = f_{n-2}/f_{n-1}$  and $N_{1}^{n-2}\cdot 1/2 = f_{n}/f_{n+1}$.   Thus,  if $n\ge2$ is even,  we have for $n\ge2$, 
\[\mathcal T^{n-1} \mathcal A_n = [0,1/2) \times \begin{cases} (f_{n-2}/f_{n-1}, f_n/f_{n+1}]\;\; &\mbox{if}\;n\; \mbox{is even};\\
\\
             ( f_n/f_{n+1},f_{n-2}/f_{n-1}]\;\; &\mbox{otherwise}.
   \end{cases}
   \]

 For $n >1$, the sets $\mathcal T^{n-1} \mathcal A_n$  are disjoint, and their union  is 
$[0,1/2) \times (0,1]\setminus \{ \frac{\sqrt{5}-1}{2}\}$; since $\mathcal A_1 = [1/2, 1] \times [0,1]$,    the union of the $\mathcal T^{n-1} \mathcal A_n$ for $n\ge 1$ is indeed $\Omega$, up to measure zero.   

With this claim, and the fact that $\mathcal G$ is given by applying $\mathcal T$ to $\mathcal T^{n-1} \Omega$, we have that $\mathcal G$ is bijective, up to measure zero.   Since $\mu$ is an invariant measure for $\mathcal T$,  it is also an invariant measure for $\mathcal G$.   Since $\mathcal T$ on $\Omega$ is a natural extension for $T$,   the images under the various $\mathcal T^n$, $n \in \mathbb N$ of  the pull-back of the Borel sigma-algebra on $I$ gives the Borel sigma-algebra on $\Omega$.   Recall (see the proof of Theorem~1 of \cite{KraaikampSSteiner:12} on p. 2219 there) that this holds since points can be separated by the integral powers of $\mathcal T$; it is easily show that this is also true for $\mathcal G$, and from this it follows that $\mathcal G: \Omega \to \Omega$ does give the natural extension of $g:T \to I$.  

Finally,  we have that $\Omega$ is double covered by $\widetilde \Sigma$ in such a way that on fibers above  $\mathcal A_n$ the $n^{\text{th}}$-return by the geodesic flow projects to give $\mathcal G$.  
     
\end{proof}

\section{Hecke triangle surfaces and Rosen fractions}\label{s:RosenCF}    An analysis of the proofs above (combined with using the aforementioned result of Sullivan \cite{Sullivan:84} in fuller generality --- see the Theorem on p.~276 there),  gives the result of this section.      

\subsection{Rosen fractions are factors of cross sections}   The Rosen continued fractions \cite{Rosen54} are defined as follows . 
Let $q\in \Z, q \geq 3$ and $\lambda = \lambda_{q}= 2 \cos \frac{\pi}{q}$. Put $\mathbb I_{q} : = [\, -\lambda/2,\lambda/2\,)\,$ and define the map $T_{q}:\mathbb I_{q} \to \mathbb I_{q}$ by
\begin{equation}
\label{def: rosenT}
T_{q}(x) := \left| \frac{1}{x} \right| - \lambda \left \lfloor\,  \left| \frac{1}{\lambda x} \right| + 1/2\right\rfloor, \textrm{ for } x\neq 0; \, T_{q}(0):=0.
\end{equation}
For  $x \in\mathbb I_{q}\,$,  put $d(x) := d_{q}(x) = \left\lfloor \left| \frac{1}{\lambda x} \right| + 1/2 \right\rfloor$ and  as usual,
$\varepsilon(x) :=\rm{sgn}(x)$. Furthermore, for $n \geq 1$ with $T_{\alpha}^{n-1}(x) \neq 0$ put
\[
\varepsilon_{n}(x)=\varepsilon_n = \varepsilon(T^{n-1}_\alpha (x)) \textrm{ and } d_n(x)=d_n=d(T^{n-1}_\alpha (x)).
\]
This yields the {\em $\alpha$-Rosen continued fraction}  of $x\,$:  
\[
x = \displaystyle{\frac{\varepsilon_1}{d_1\lambda +  \displaystyle{\frac{\varepsilon_2}{d_2\lambda + \dots}}}} 
\,,
\]
where $\varepsilon \in \{ \pm 1\}$ and $d_i \in \N$.    
\\
Rosen introduced his continued fractions to study the Hecke  groups.   The Hecke (triangle Fuchsian) group $G_q$ with $q \in \{3, 4, 5, \dots\,  \}$ is the group generated by 
\[
 \begin{pmatrix}
	1   &  \lambda_q\\
 	0   &  1
\end{pmatrix}\, \, \text{and} \;\; 
 \begin{pmatrix}
        0   &  -1\\
        1   &  0
\end{pmatrix},\]
with $\lambda_q$ as above.

\begin{Thm} \label{t:rosenCF}   For each $q \ge 3$,  the Rosen fraction map of index $q$ is given as the factor of a cross section of the geodesic flow on the unit tangent bundle of $G_q\backslash \mathbb H$.
\end{Thm}

\begin{proof}     Due to the general nature of Sullivan's result,  Equation~\eqref{e:entropyReturn} holds for the first return map to any cross section $\mathcal C$ for the geodesic flow on unit tangent bundle of $G_q\backslash \mathbb H$ for any $q$.     For each $q$,  Burton-Kraaikamp-Schmidt \cite{BurtonKraaikampSchmidt00} determined a planar natural extension on a region $\Omega_q$ with the measure $\mu$ as above.   It is easily verified that the analog of Equation \eqref{e:biject} holds, and we thus
define $\widetilde \Sigma_q$ completely analogously to the $\widetilde \Sigma_{\alpha}$ above.     The analog of  Lemma~\ref{l:distinctClasses} goes through,  as for any fixed $q$,  all elements of $G_q$ have their entries lying in the (countable!) algebraic number field $\mathbb Q(\lambda_q)$.      The analog of Lemma~\ref{l:flowFormula} holds upon replacing each occurrence of $d$ there by the appropriate $d$ times  $\lambda_q$.      Finally,  the result holds,  since Nakada \cite{Nakada10} showed that the entropy of the Rosen continued fraction map equals  one-half times the quotient of the volume of the unit tangent space of $G_q\backslash \mathbb H$ by the $\mu$-area of $\Omega_q$.  
\end{proof}

\subsection{Cross sections and $\alpha$-Rosen fractions} 
Dajani, Kraaikamp and Steiner \cite{DajaniKraaikampSteiner09} introduced  the $\alpha$-Rosen fractions, a generalization combining the idea of Nakada's $\alpha$-continued fractions with Rosen's continued fractions.   

Let $q\in \Z, q \geq 3$ and $\lambda = \lambda_{q}= 2 \cos \frac{\pi}{q}$. For $\alpha \in \left[\,0, \frac{1}{\lambda}\,\right],$ we put $\mathbb I_{q,\alpha} : = [\, \lambda(\alpha-1),\lambda\alpha\,)\,$ and define the map $T_{\alpha}:\mathbb I_{q,\alpha} \to \mathbb I_{q,\alpha}$ by
\begin{equation}
\label{def: Talpha}
T_{\alpha}(x) := \left| \frac{1}{x} \right| - \lambda \left \lfloor\,  \left| \frac{1}{\lambda x} \right| + 1 -\alpha \right\rfloor, \textrm{ for } x\neq 0; \, T_{\alpha}(0):=0.
\end{equation}
For  $x \in\mathbb I_{q,\alpha}\,$,  put $d(x) := d_{\alpha}(x) = \left\lfloor \left| \frac{1}{\lambda x} \right| + 1 -\alpha \right\rfloor$ and  
$\varepsilon(x) :=\rm{sgn}(x)$. Furthermore, for $n \geq 1$ with $T_{\alpha}^{n-1}(x) \neq 0$ put
\[
\varepsilon_{n}(x)=\varepsilon_n = \varepsilon(T^{n-1}_\alpha (x)) \textrm{ and } d_n(x)=d_n=d(T^{n-1}_\alpha (x)).
\]
The  {\em $\alpha$-Rosen continued fraction}  of $x$ is then defined in what now is the obvious fashion.     Fixing $\alpha = \frac{1}{2}\,$  results in the Rosen fractions.   On the other hand, fixing  $q=3\,$ and considering general $\alpha\,$, we have Nakada's $\alpha$-expansions. 

Using direct methods, similar to those of \cite{BurtonKraaikampSchmidt00}  for the classical Rosen fractions, planar natural extensions $\Omega_{q, \alpha}$ with invariant measure $\mu$ as above for certain of the $\alpha$-Rosen fractions are given in  \cite{DajaniKraaikampSteiner09}.   For each index $q$, Kraaikamp-Schmidt-Smeets \cite{KraaikampSSmeets10} determine the value $\alpha_0 = \alpha_0(q)$ such that $[\alpha_0, 1/\lambda]$ is the maximal interval containing $1/2$ with  $\Omega_{q,\alpha}$ being connected for each value of $\alpha$ in this interval.     The ``quilting'' of $\Omega_{q, 1/2} = \Omega_q$ to determine these $\Omega_{q, \alpha}$ is then fairly straightforward;   \cite{KraaikampSSmeets10} find a subinterval containing $1/2$ on which entropy is constant,  since also $\mu$-measure of $\Omega_{q, \alpha}$ is constant.   They also give an argument (see Lemma~12 there) that easily implies the analog of Equation~\eqref{e:entropArea} for these values of $q$ and $\alpha$.    Combining these results with Nakada's entropy calculation for the classical Rosen maps, we have the following.    

\begin{Thm} \label{t:rosenAlpCF}   For each $q > 3$,  let $\alpha_0(q)$ be as in \cite{KraaikampSSmeets10}.   Then for each $q$ and  $\alpha \in 
[\,\alpha_0(q), 1/\lambda_q\,]$  the $\alpha$-Rosen fraction map of index $q$ is given as the factor of a cross section of the geodesic flow on the unit tangent bundle of $G_q\backslash \mathbb H$.
\end{Thm}

\section{Additional remarks}

\subsection{The volume of the unit tangent bundle on the modular surface}

This volume, for the normalization of the Haar measure used in this paper, can be derived by an elementary computation.

Indeed, when we consider the special case of $T_1$, the classical Gauss map, we see that one can find a cross-section of the geodesic flow which is a double cover of a surface parametrized by $\{(x,y)|0\le x\le 1, 0\le y\le \frac 1{1+x}\}$. This cross-section defines a fundamental domain of the action of $\text{PSL}_2(\mathbb Z)$ on $\text{PSL}_2(\mathbb R)$. 

With this parametrization, the return time has been seen above to be $-2\log x$, and the invariant measure is Lebesgue measure; hence the measure of the unit tangent bundle is $2\int_0^1\int_0^{\frac 1{1+x}}-2\log x\, dy\, dx=\frac{\pi^2}3$.

\subsection{Variants of continued fractions with determinant 1}
\begin{figure}[h]
\includegraphics[width=5cm]{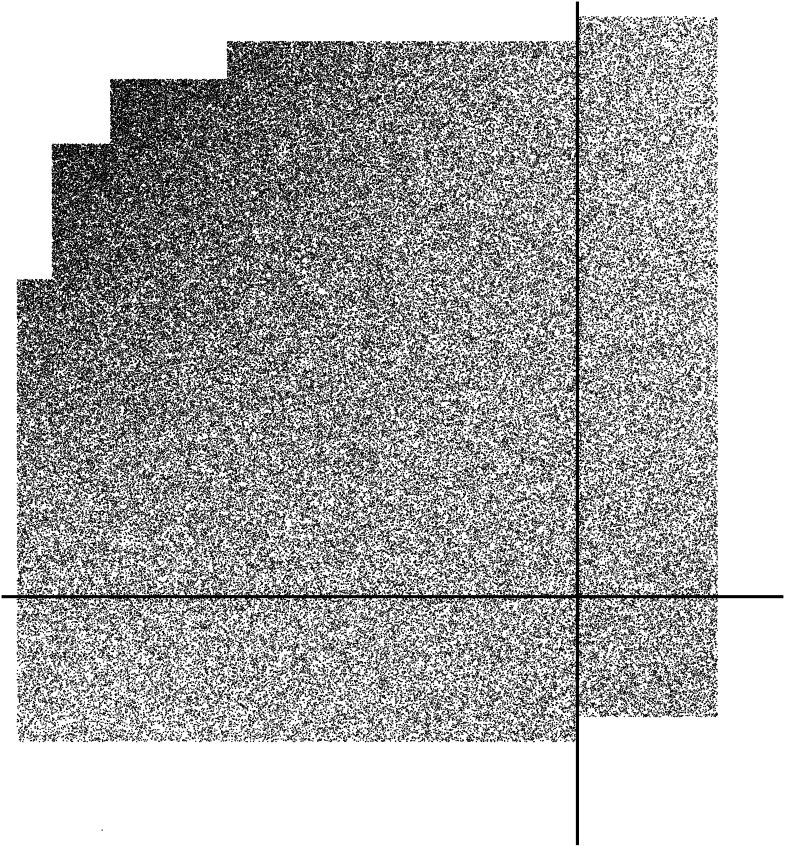}
\includegraphics[width=5cm]{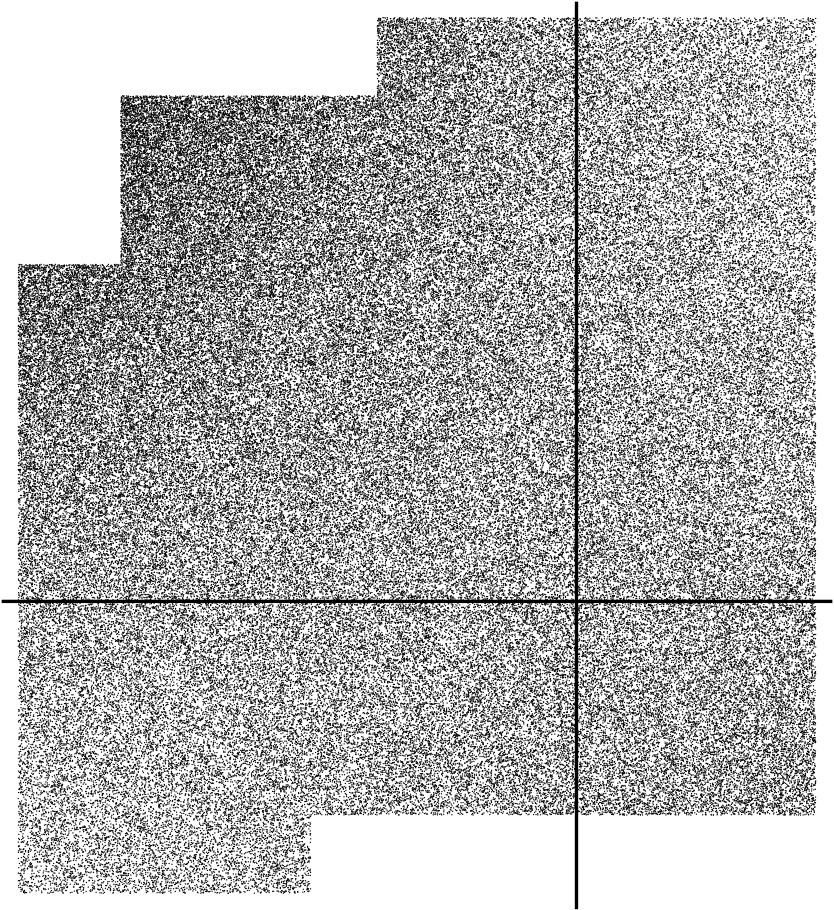}
\includegraphics[width=5cm]{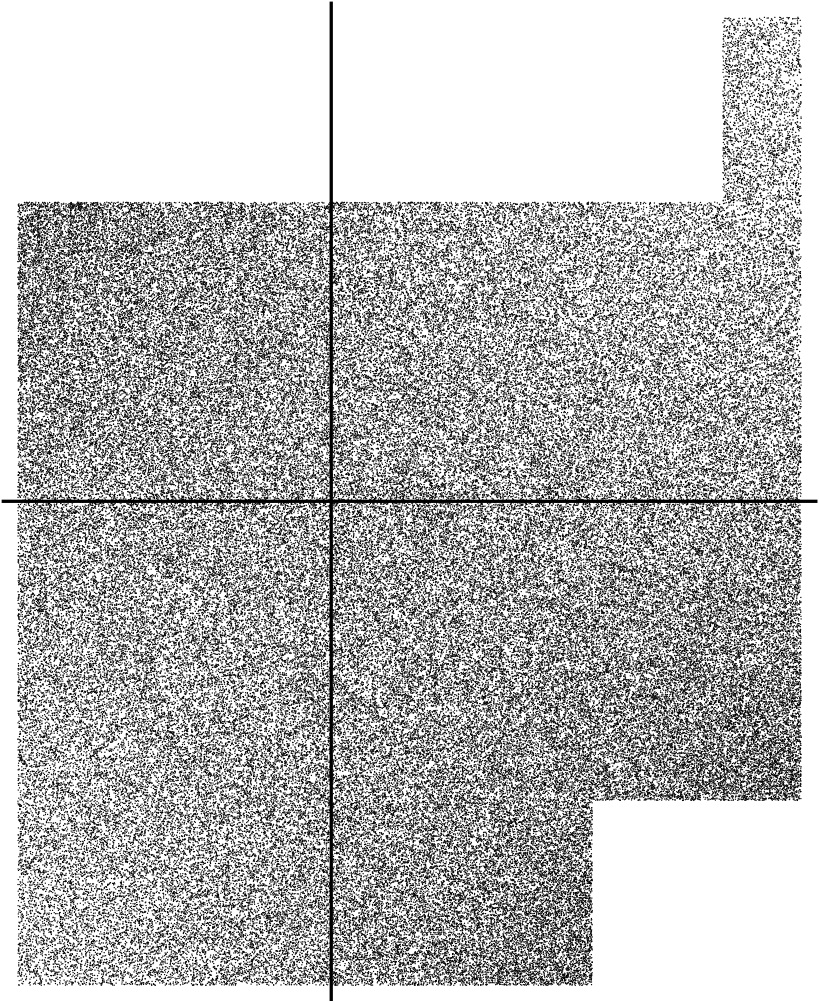}
		\caption{The domain of the natural extension of the continued $\alpha$-fraction with positive determinant, for $\alpha=0.2, \, 0.3, \,0.6$}
	\label{fig:dessExtp}
\end{figure}

A technical difficulty  in all this work arises from the fact that the linear maps underlying the transformations $T_{\alpha}$ are in $\text{GL}_2(\mathbb Z)$ and not necessarily in  $ \text{SL}_2(\mathbb Z)$;  they can have determinant -1. Thus the necessity of the 2-fold covering of the natural extension ---   since $ \text{GL}_2(\mathbb Z)$ does not act on the hyperbolic plane, we alternate between the two sheets of the covering when the determinant is -1. The only case when this problem does not occur is for $T_0$.

\begin{figure}[h]
\includegraphics[width=5cm]{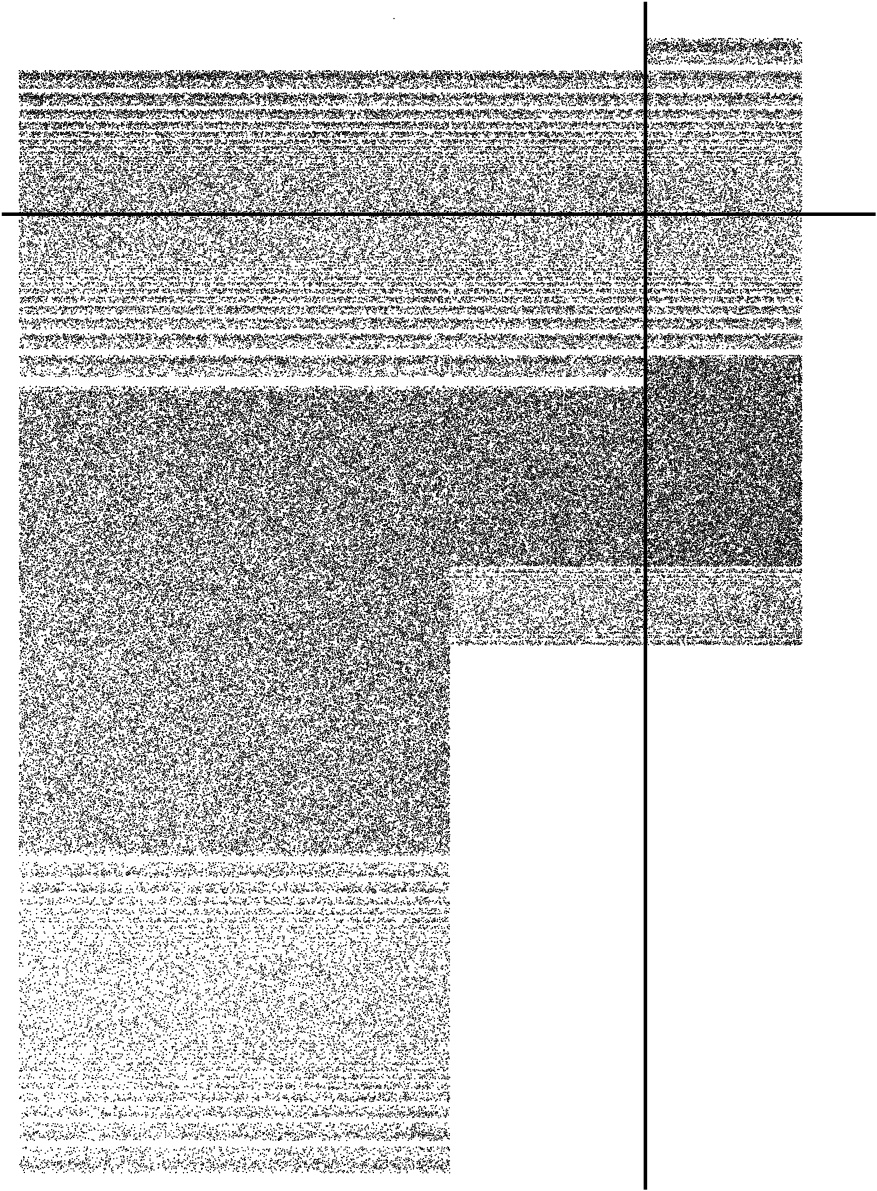}
\includegraphics[width=5cm]{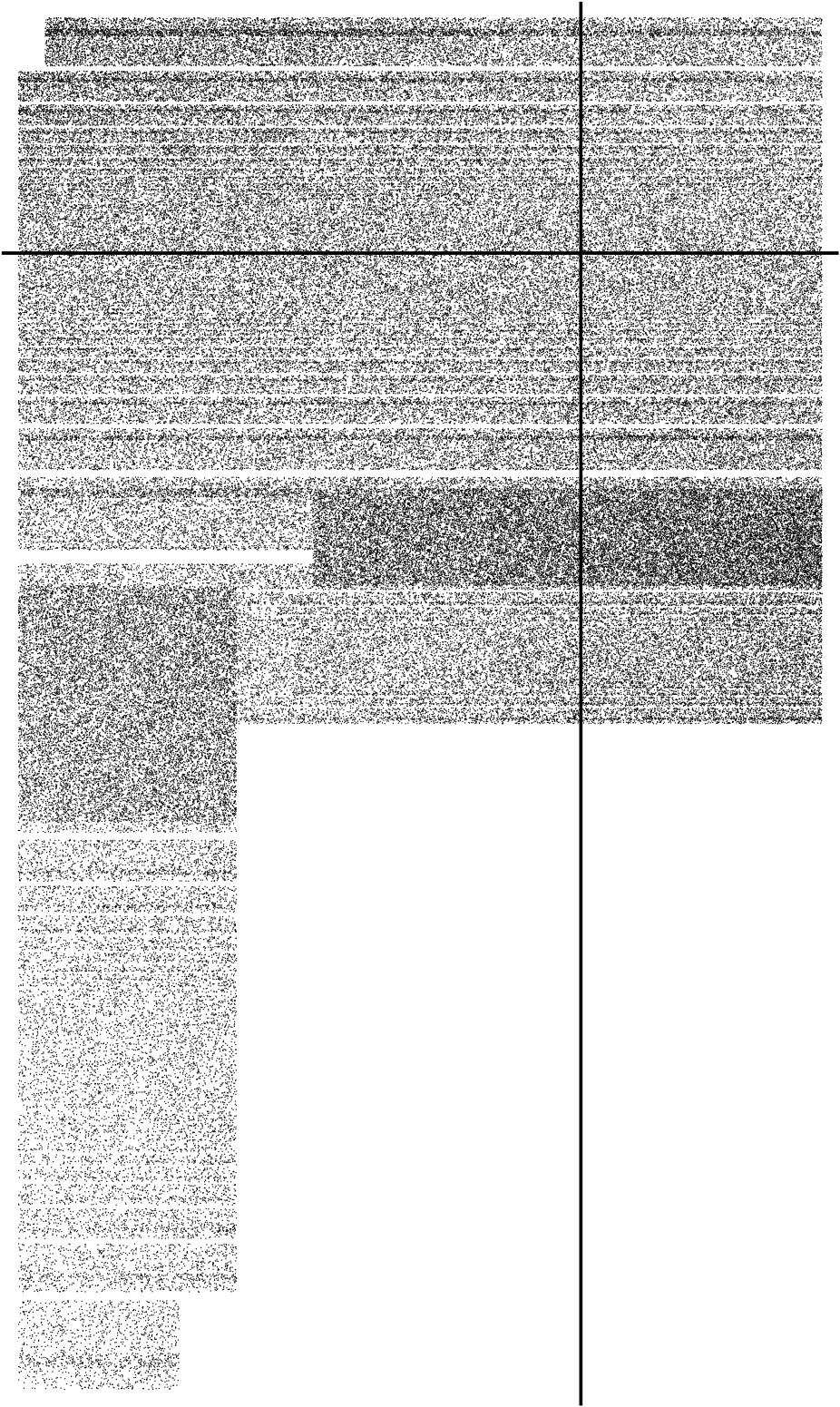}
\includegraphics[width=5cm]{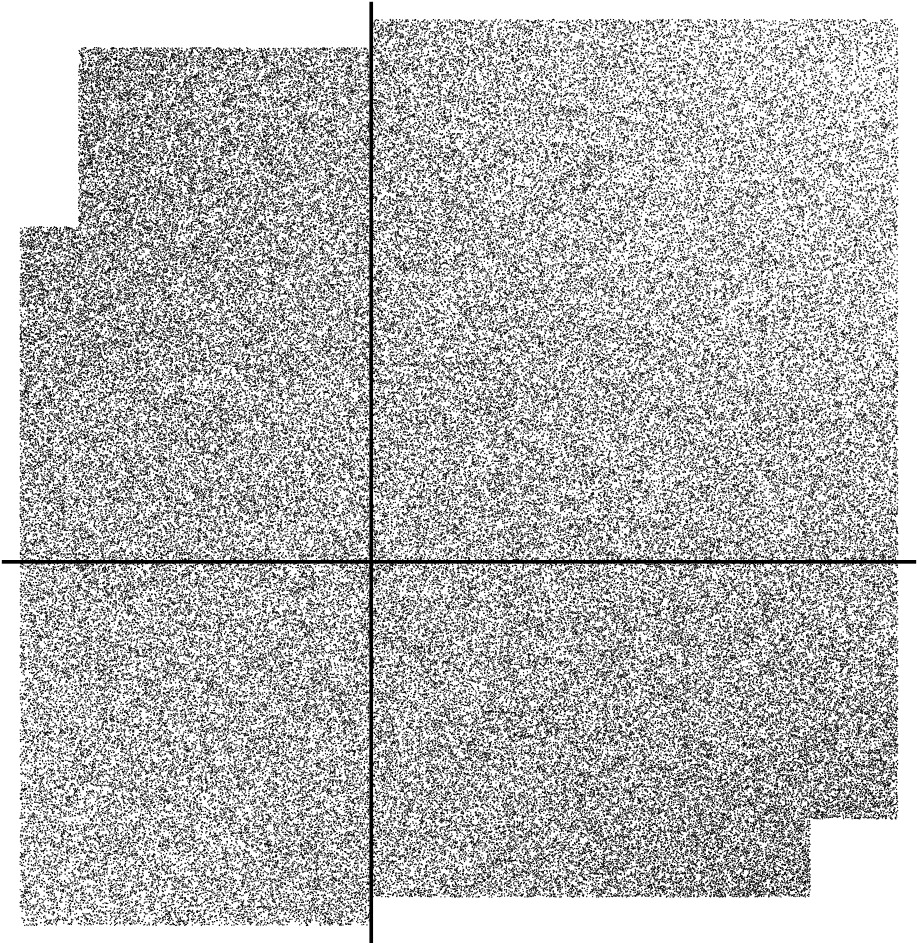}
		\caption{The domain of the natural extension of the continued $\alpha$-fraction with  negative determinant, for $\alpha=0.2, \, 0.3, \,0.6$}
	\label{fig:dessExtm}
\end{figure}
One might  avoid this problem by defining a variant of the $\alpha$-continued fraction as:
  $$S_{\alpha}(x) = \frac{-1}{x}  -  \left \lfloor\,  \frac{-1}{x}  + 1 -\alpha \right\rfloor \text{ for (non-zero) }x \in  [1-\alpha, \alpha)$$

In that case, the fractional linear maps associated with the transformations are always in   $ \text{SL}_2(\mathbb Z)$, so we will not need a double cover of the natural extension. The limit cases $S_0$ and $S_1$ both have an indifferent fixed point, which is probably the reason why this family has not been very much studied (although $S_1$ is known as the \emph{backwards continued fraction}).

Experimental studies indicate that the natural extension behaves much more nicely than for the $\alpha$-continued fraction; see Fig.~\ref{fig:dessExtp}, where we have shown the first 200 000 iterates of the point $(e/10, 0)$. In particular, its domains seems to be connected for all $\alpha\in (0,1)$, and any vertical line intersects it in an interval, the domain being bounded by two increasing step functions.   Hence its study might be easier than that of the $\alpha$-continued fraction.     (We note that Nakada's ~\cite{Nakada92}  variant of Rosen continued fractions,  mentioned in our introduction,  also has all associated fractional linear maps of  determinant one.)

By contrast,  the functions $ \frac{1}{x}  -  \left \lfloor\,  \frac{1}{x}  + 1 -\alpha \right\rfloor \text{ for (non-zero) }x \in  [1-\alpha, \alpha)$, which are perhaps more natural, seem to be more delicate to study in the general case, see Fig. ~\ref{fig:dessExtm}, where we show the first 200 000 iterates of the same point; the corresponding fractional linear maps always have determinant -1, and for small values of $\alpha$, the domain of the natural extension appears to be disconnected, the vertical sections being Cantor-like with nonempty interior.

\end{document}